\documentclass[12pt]{amsproc}    
\usepackage{latexsym,amsmath,amssymb,amsfonts,amscd,amsthm,amstext} 
\usepackage{tikz}

\swapnumbers
           \theoremstyle{definition}
\newtheorem{defC}{Definition}[section]
\newtheorem{exer}[defC]{Exercise}

           \theoremstyle{plain}
\newtheorem{lemma}[defC]{Lemma}
\newtheorem{thm}[defC]{Theorem}    
\newtheorem{cor}[defC]{Corollary}                                                   
\newtheorem{conjecture}[defC]{Conjecture}

           \theoremstyle{remark}
\newtheorem{com}[defC]{Comment}

\begin{document}
\title[Monoidal category of graphs]{Monoidal category of operad of graphs}
\author[Ch\'avez Rodr\'{\i}guez]{Mar\'{\i}a Ernestina Ch\'avez Rodr\'{\i}guez}
\address{Universidad Nacional Aut\'onoma de M\'exico, Facultad de Estudios Superiores,  C.P. 54714 Cuautitl\'an Izcalli, Apartado Postal \# 25, Estado de M\'exico}
\email{mechavez@servidor.unam.mx}
\author[Oziewicz]{Zbigniew Oziewicz}
\address{Department of Mathematics, The University of Texas at San Antonio, One UTSA Circle, San Antonio, Texas 78249-0624}
\curraddr{Universidad Nacional Aut\'onoma de M\'exico, Facultad de
Estudios Superiores,  C.P. 54714 Cuautitl\'an Izcalli, Apartado Postal
\# 25, Estado de M\'exico}\email{oziewicz@unam.mx}
\date{October 17, 2010} 
\thanks{Proceedings of the Conference on Mathematical Sciences for Advancement of Science and Technology, MSAST, Kalkata (Calcutta), India, December 2010. Institute of Mathematics, Bio-informatics, Information Technology and Computer-science IMBIC, www.imbic.org/index.html}
\thanks{This work is Supported by Programa de Apoyo a Proyectos de Investigaci\'on e Innovaci\'on Te\-cno\-l\'ogica, UNAM, Grant PAPIIT \# IN104908, 2008--2010}

\newcommand{\N}{\mathbb{N}}\newcommand{\C}{\mathbb{C}}
\newcommand{\R}{\mathbb{R}}\newcommand{\Z}{\mathbb{Z}}
\newcommand{\cM}{\mathcal{M}}
\newcommand{\cA}{\mathcal{A}}\newcommand{\cC}{\mathcal{C}}
\newcommand{\cF}{\mathcal{F}}
\newcommand{\ie}{\textit{i.e.}$\,$}
\newcommand{\be}{\begin{equation}}
\newcommand{\ra}{\longrightarrow}
\newcommand{\id}{\operatorname{id}}
\newcommand{\cat}{\operatorname{\text{\textbf{cat}}}}
\newcommand{\bimod}{\operatorname{\text{-\textbf{bimod}}}}
\newcommand{\obj}{\operatorname{obj}}
\newcommand{\lin}{\operatorname{lin}}
\newcommand{\emod}{\operatorname{mod}}
\newcommand{\alg}{\operatorname{alg}}
\newcommand{\der}{\operatorname{der}}
\newcommand{\cog}{\text{cog}}
\newcommand{\aut}{\text{aut}}
\newcommand{\End}{\operatorname{End}}
\newcommand{\Nat}{\operatorname{Nat}}
\newcommand{\sh}{\text{sh}}
\newcommand{\hopf}{\text{hopf}}
\newcommand{\Span}{\text{span}}
\newcommand{\im}{\operatorname{im}}
\newcommand{\gen}{\operatorname{gen}}
\newcommand{\grade}{\operatorname{grade}}
\newcommand{\Cl}{\mathcal{C}\ell}
\newcommand{\half}{\textstyle{\frac{1}{2}}}
\newcommand{\bt}{\begin{tabular}{c}}
\newcommand{\et}{\end{tabular}}

\begin{abstract} Usually a name of the category is inherited from the name of objects. However more relevant for a category of objects and morphisms is an algebra of morphisms. Therefore we prefer to say a category of \textit{graphs} if every morphism is a graph. 

In a monoidal category every morphism can be seen as a graph, and a partial algebra of morphisms possesses a structure of an operad, operad of graphs. We consider a monoidal category of operad of graphs with underlying graphical calculus. If, in particular, there is a single generating objects, then each morphism is a bi-arity graph. The graphical calculus, multi-grafting of morphisms, is developed \textit{ab ovo}.      

We interpret algebraic logic and predicate calculus within a monoidal category of operad of graphs, and this leads to the graphical logic.

A logic based on a braided monoidal category is said to be the
braided logic. We consider a braided monoidal category generated by
one object. We are demonstrating how the braided logic is related to
implicative algebra and to the Heyting algebra (in contrary to the
Boolean algebra) and therefore must be more related to the quasigroups
then to the lattices.

Some applications to classical logic, to modal logic and to {\L}u\-ka\-sie\-wicz three-valued logic are considered.\end{abstract}\maketitle
\newpage\tableofcontents

\section{Operad of graphs} We work within a monoidal category (almost a tensor category) generated by a single object. In this case the set of all objects of the single generated  monoidal category coincides with the set of non-negative integers (with the set of natural numbers) $\N.$ Thus the set of all objects is, $\obj\cat=\N,$ and $1\in\N$ is a generating object.

Therefore each morphism (an arrow of a category) is characterized by a \textit{pair} of non-negative integers, morphims are bi-graded, and we refer to this pair  \{input, output\} = \{entrance, exit\}, as to the type or arity of the morphism = \{arity-in, arity-out\},
\begin{gather}\N\ni m\quad\xrightarrow{\quad\text{morphisms of $(m\rightarrow n)$-arity\quad}}\quad n\in\N.\end{gather}

Garrett Birkhoff in his \textit{Lattice Theory} in 1940, within Universal algebra considered an $n$-ary operation symbol to be carefully distinguished from a model of the operation. 
Traditionaly, every operation symbol in an algebra is
of arity (of type) $\in\N,$ \ie of arity $\N\ni n\mapsto 1$ with one exit. 
Every such operation is considered here to be a morphism in a monoidal category. For example the traditional $0$-ary (or null-ary) operation  in our terminology must be bi-graded. It is a $(0\rightarrow 1)$-operation/morphism, $\in\cat(0,1)$ Similarly traditional $n$-ary operations are seen here as morphisms $\in\cat(n,1).$

We need to introduce the graphical notation we are using. Every morphism $\in\cat(m,n),$ $m,n\in\N,$ is visualized as a node with a number of outer leaves (representing the source and target objects): on top the input object of $m$-arity-in, and at bottom the output object of $n$-arity-out as illustrated in Figure \ref{morphism}. Throughout this paper the graphs are directed and we read them from the top to the bottom. Bi-gradation is in accordance with the time arrow. Reading these graphs upside down, from the bottom to the top, can be
equivalent to the inverse of the time's arrow.

For example, a morphism $\in\cat(2,1)$ in Figure \ref{morphism} is traditionally called a binary algebra. Each morphism is a graph with no outer nodes. Every node is inner, and every edge is an outer leave or an inner edge, which represents an object of a monoidal category. 

\begin{figure}[h!]\begin{center}
\begin{tikzpicture}[line width=1pt]
\draw(-2,0) .. controls (-1.7,-1) and (-1.3,-1) .. (-1,0);\draw(-1.5,-0.8)--(-1.5,-1.4);
\draw(-1.5,0.4)node{$\overbrace{\hspace{3mm}\text{input: $2$}\hspace{3mm}}$};
\draw(2.5,0)--(3.5,-0.8);\draw(3,0)--(3.5,-0.8);
\draw(4,0)node{$\ldots$};
\draw(4.5,0)--(3.5,-0.8);
\draw(3.5,0.4)node{$\overbrace{\hspace{7mm}\text{input: $m$}\hspace{7mm}}$};
\draw(3.5,-0.8)--(2.5,-1.6);
\draw(3.5,-0.8)--(3,-1.6);\draw(3.5,-0.8)--(4.5,-1.6);\draw(4,-1.6)node{$\ldots$};
\draw(3.5,-2)node{$\underbrace{\hspace{7mm}\text{output: $n$}\hspace{7mm}}$};
\end{tikzpicture}\end{center}
\caption{A $(2\rightarrow 1)$-morphism $\in\cat(2,1),$ and $(m\rightarrow n)$-morphism $\in\cat(m,n)$ is an arrow from $m$ to $n.$\label{morphism}}\end{figure}
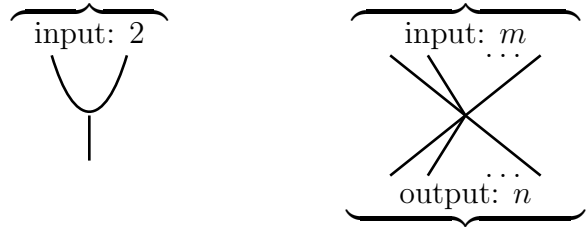

Two morphisms with the same arity-in and arity-out are said to be \textit{paralell}, they are parallel arrows, $m\rightrightarrows n.$ 

\begin{defC}[Schizophrenic] In monoidal category $\cat,$ we postulate a set, $\Omega\equiv\cat(0,0),$ to be a candidate for a dualizing or schizophrenic set. A plant $0\mapsto 0$ can be always blot out.\end{defC}

Morphisms can be composed in many different ways, like garden plants are grafted. Composition of morphisms generalizes the construction of words from an alphabet, however words are constructed by concatenations only, whereas our bi-graded morphisms allow diverse `non-linear' compositions.  

We refer to each composition as a grafting. The concatenation or  juxtaposition of arrows is considered also as a special grafting. For example a concatenation of an arrow from $\cat(2,1)$ with another arrow from $\cat(2,1),$ gives an arrow from $\cat(4,2),$ see Figure \ref{concatenation}.
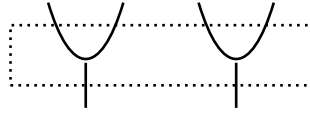
\begin{figure}[h!]\begin{center}
\begin{tikzpicture}[line width=1pt][line width=1pt]
\draw(-2,0) .. controls (-1.7,-1) and (-1.3,-1) .. (-1,0);\draw(-1.5,-0.8)--(-1.5,-1.4);
\draw(0,0) .. controls (0.3,-1) and (0.7,-1) .. (1,0);\draw(0.5,-0.8)--(0.5,-1.4);
\draw[dotted](-2.5,-1.1)rectangle(1.5,-0.3);
\end{tikzpicture}\end{center}
\caption{The concatenation of $Y\in\cat(2,1)$ with $Y$ gives an arrow in $\cat(4,2).$\label{concatenation}}\end{figure}  

Occasionally a bigraded morphism we will call in a diverse way, a graph, a garden \textit{plant} with root and crown of branches, operation, symbol, letter, atom, a formal predicate. 
The multitude of names is better expressing what is going on. 

Except for concatenations that will be not used in what follows, each grafting is also bigraded, and $(i\rightarrow j)$-graft means joining of the $i$-th output of the first $(m\rightarrow n\neq 0)$-morphism, $i\in\{1,2,\ldots,n\},$ with $j$-th input of the second $(0\neq k\rightarrow l)$-morphism, $j\in\{1,2,\ldots,k\},$ as illustrated in several examples in next Figures.  

`The plants' can be grafted in more ways then `the letters' can be juxtapositioned.
Reader can freely choose the favorite name.
 
\begin{defC}[Free operad] Every chosen finite set of such bi-graded morphisms we treat as an alphabet generating, by all possible graftings, a free operad (clon, abstract algebra, variety). Every generator of an operad is of the type in: entrance $\N\ni m\mapsto n\in\N$ exit,
\ie all plants are bi-graded as it is the case in the nature.\end{defC}

A null-ary operation, a creator, is of arity $(0\mapsto 1)\simeq\cat(0,1).$ The null-ary
co-operation, a killer, is of type $(1\mapsto 0)\simeq\cat(1,0),$ and is unique, so that a
category possess the terminal object. The killer is the process of forgeting all, the annihilator.  An un-ary cooperation or operation is of type $(1\mapsto 1)\simeq\cat(1,1),$ and the general type is: entrance $\N\ni m\mapsto n\in\N$ exit.

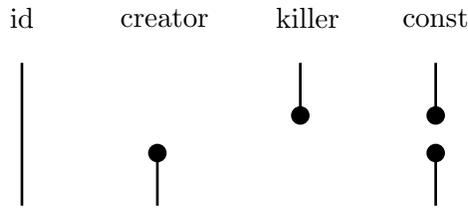
\begin{figure}[!ht]\begin{centering}\small
\begin{tikzpicture}[line width=1pt][line width=1pt]
\draw(-1,0.6)--(-1,-1.3);  \draw (-1,1.2) node{$\id$};
   \draw(.9,1.2) node{creator};
   \draw( .8,-1.3)-- ( .8,-.6);\filldraw ( .8,-.6) circle (3pt);
 
   \draw (2.8,1.2) node{killer};
   \draw (2.7,-.1) -- (2.7,.6);
   \filldraw (2.7,-.1) circle (3pt);
 
   \draw (4.5,1.2) node{const};
   \draw (4.5,-.1) -- (4.5,.6);
   \filldraw (4.5,-.1) circle (3pt);
   \draw (4.5,-1.3) -- (4.5,-.6);
   \filldraw (4.5,-.6)circle (3pt);
\end{tikzpicture} \end{centering}
 \caption{From the left: $\id\simeq\arrowvert\in\cat(1,1),$ creator $\in\cat(0,1),$ forgetting $\in\cat(1,0),$ and composed the constant morphism $\in\cat(1,1).$}\label{music}\end{figure}
  
\begin{figure}[!ht]\begin{centering}\small
   \begin{tikzpicture}[line width=1pt][line width=1pt]\draw (0.6,1.2) node{$nil$potent};
   \draw (0,-1.3) -- (0,.6);
   \filldraw (0,-.1) circle (3pt);
   \filldraw (0,-.6) circle (3pt);
   \draw (0.6,-.4) node{$\sim$};
   \draw (1.2,-.1) --(1.2,.6);
   \filldraw (1.2,-.1) circle (3pt);
   \draw (1.2,-1.3) -- (1.2,-.6);
   \filldraw (1.2,-.6) circle (3pt);
 
   \draw (3.6,1.2) node{$uni$potent};
   \draw (3,-1.3) -- (3,.6);
   \filldraw (3,-.1) circle (3pt);
   \filldraw (3,-.6) circle (3pt);
   \draw (3.6,-.4) node{$\sim$};
   \draw (4.2,-1.3) -- (4.2,.6);

   \draw (6.6,1.2) node{$idem$potent};
   \draw (6,-1.3) -- (6,.6);
   \filldraw (6,-.1) circle (3pt);
   \filldraw (6,-.6) circle (3pt);
   \draw (6.6,-.4) node{$\sim$};
   \draw (7.2,-1.3) -- (7.2,.6);
   \filldraw (7.2,-.4) circle (3pt);\end{tikzpicture} \end{centering}
 \caption{Three examples of unaries $\in\cat(1,1)$ defined by relations among parallel morphisms.}\label{music2}\end{figure}
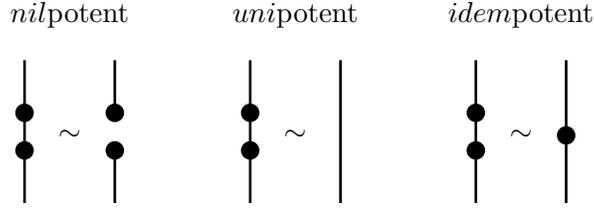

\begin{defC}[Quotient operad] An equivalence relations among parallel morphisms, among parallel arrows, determine two-sided ideal in the free operad. This leads to the quotient operad. On Figure \ref{music2}, three identities (the equivalence relations in an operad among parallel morphisms) on the right are, by definition, the abstract minimum polynomials.\end{defC}

The boolean negation is an example of the unipotent unary
operation, but nil- and idem-potents are not classical, \ie does not
exist for the two elements object $\{$true,false$\}.$

Joyal \& Street [1991] introduced the \textit{valuation} of a graph as
the pair of aplications [Joyal\& Street 1991, Definition 1.3, p. 64],
\begin{gather*}\text{edges}\ra\obj\cat,\qquad\text{nodes}\ra
\text{arrows}\cat.\end{gather*}

In our convention graph has no outer nodes. Every node is inner and
represent a functor or a map, \ie a process, an action, a co-action,
operation, multi-functor, evaluation, function, etc. Every edge,
including outer leaves, represent object.

A bifunctor of bin-ary operation, an anihilation $\in(2\mapsto 1)\simeq\cat(2,1),$
denoted by two initial leaves and one node. A binary co-operation is a
decomposition (of information), a splitting, duplication with the
mutants, procreation process $\in(1\rightarrow 2)\simeq\cat(1,2).$ The bin-aries:
anihilation $\in(2\mapsto 1)\simeq\cat(2,1),$ creation $\in(1\mapsto 2)\simeq\cat(1,2),$ and
scattering $\in(2\mapsto 2)\simeq\cat(2,2),$ are represented by the prime graph
nodes on Figure \ref{binope}.

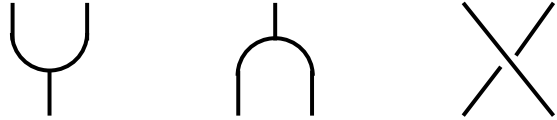
\begin{figure}[!ht]
 \begin{centering}
   \begin{tikzpicture}[line width=1.5pt]
     \draw (.5,-1.3) -- (.5,-1.9);
     \draw (0.009,-.9) -- (.009,-.4);
     \draw (1,-.9) -- (1,-.4);
     \draw[ ] (0mm,-8mm) arc (180:360:5mm);
    \draw (3.5,-.4) -- (3.5,-.9);
     \draw (3.009,-1.3) -- (3.009,-1.9);
     \draw (3.99,-1.3) -- (3.99,-1.9);
     \draw (30mm,-13.7mm) arc (180:0:5mm);
\draw (6.,-.4) -- (7.2,-1.9);
     \draw (7.2,-.4) -- (6.7,-1.1);
     \draw (6,-1.9) -- (6.5,-1.25);
   \end{tikzpicture}
  \end{centering}
  \caption{The binaries: operation (anihilation), co-product (procreation) and a
   scattering (prebraid).}
  \label{binope}
 \end{figure}

\begin{defC}[The plication] A pro-creation process $\in(1\mapsto n\geq
2)\simeq\cat(1,n),$ in realization the process of the copy of an identical variables,
informations, thought, ideas, things, species, genuses, $\ldots$ [Reader
can freely choose the favorite name], is called $n$-plication
(du-plication, tri-plication, multi-plication) or $n$-ary plication.
The (du)-plication is co-associative and co-commutative, is the
`group-like' co-operation like reproduction in biology, like mitosis,
cariocinesis $\triangle:a\mapsto(a,a).$ Compare with [Lambek \& Scott,
page 8, Exercise 3]. Abstractly (two-sided) duplication is defined by
the pair of relations on Figure \ref{mitosis}.

\begin{figure}[!ht]
 \begin{centering}
   \begin{tikzpicture}[line width=1pt]  
     \draw (-1.5,-.4) -- (-1.5,-.9);
     \draw (-2.009,-1.3) -- (-2.009,-1.9);
     \draw (-0.99,-1.3) -- (-0.99,-1.7);
     \draw (-20mm,-13.7mm) arc (180:0:5mm);
     \filldraw (-0.99,-1.6) circle (3pt);
 \draw (2.,-.4) -- (2,-1.9);
\draw ( .5,-1.1) node(x)  {$\sim$};
     \draw (3.5,-1.1) node(x)  {$\sim$};
     \draw (5.5,-.4) -- (5.5,-.9);
     \draw (5.009,-1.3) -- (5.009,-1.7);
     \draw (5.99,-1.3) -- (5.99,-1.9);
     \draw (50mm,-13.7mm) arc (180:0:5mm);
     \filldraw (5.009,-1.6) circle (3pt);
\end{tikzpicture}\end{centering}\caption{Mitosis.}\label{mitosis}\end{figure}
\end{defC}

\begin{com} In realization on the set $S,$ mitosis is known as the
diagonal or identity relation on $S,$ mitosis $\in 2^{(S\times S)}.$
The name \textit{multiplication} we use for the co-process
$a\mapsto(a,a,\ldots,a),$ contrary to the usual meaning in literature of
the binary operation.
\end{com}

Every binary cooperation is represented diagrammatically by node with
two exit edges (two exit streams) as shown on Figures \ref{mitosis} and \ref{cooper}.

Every general co-operation $\in\cat(1,n),$ like meiosis in biology
$\in\cat(1,4),$ can be considered as the composition of the multiple
mitosis with a set of $n$ (different) modal unary operations which we
interpret as the \textit{mutant} operations for cooperations.

\begin{figure}[!ht]
 \begin{centering}
   \begin{tikzpicture}[line width=1pt]
     \draw (-0.5,-.4) -- (-0.5,-.9);
     \draw (-1.009,-1.3) -- (-1.009,-1.9);
     \draw (-0.00,-1.3) -- (-0.00,-1.9);
     \draw (-10mm,-13.7mm) arc (180:0:5mm);
     \draw (-.5,-.9) circle (8pt);
     \draw (-0.7,-.7) -- (-0.25,-1.1);
     \draw (-0.55,-.6) -- (-0.24,-0.9);
     \draw (-0.8,-.85) -- (-0.37,-1.2);

     \draw (1.8,-1.1) node(x) {$\sim$};
  
     \draw (4.5,-.4) -- (4.5,-.9);
     \draw (4.009,-1.3) -- (4.009,-1.9);
     \draw (4.99,-1.3) -- (4.99,-1.9);
     \draw (40mm,-13.7mm) arc (180:0:5mm);
     \filldraw (4.009,-1.3) circle (3pt);
     \filldraw (4.99,-1.3) circle (3pt);
\end{tikzpicture}\end{centering}
\caption{Every cooperation (read from the top) is the grafting (composition)
of the duplication with the pair of unary mutants.}\label{cooper}\end{figure}
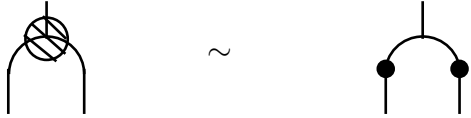

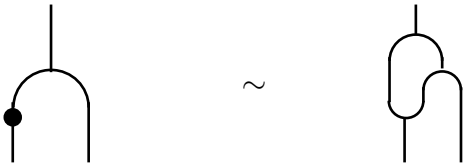
\begin{figure}[!ht]
 \begin{centering}
   \begin{tikzpicture}[line width=1pt]
     \draw (-0.5,-.0) -- (-0.5,-.9);
     \draw (-1.009,-1.3) -- (-1.009,-2.1);
     \draw (-0.00,-1.3) -- (-0.00,-2.1);
     \draw (-10mm,-13.7mm) arc (180:0:5mm);
     \filldraw (-1.009,-1.5) circle (3pt);

     \draw (2.2,-1.1) node(x) {$\sim$};
  
     \draw (4.35,-0) -- (4.35,-.4);
     \draw (4.,-.7) -- (4.,-1.3);
     \draw (4.7,-.7) -- (4.7,-.85);
     \draw (40mm,-7.5mm) arc (180:0:3.5mm);
     \draw (44.5mm,-11.4mm) arc (180:0:2.5mm);
     \draw (4.45,-1.1) -- (4.45,-1.3);
     \draw (4.95,-1.1) -- (4.95,-2.1);
     \draw (39.9mm,-12.8mm) arc (180:360:2.3mm);
     \draw (4.2,-1.5) -- (4.2,-2.1);
   \end{tikzpicture}
  \end{centering}
  \caption{A right mutant from binary operation.}
  \label{rightmut}
\end{figure}

\begin{figure}[!ht]
 \begin{centering}
   \begin{tikzpicture}[line width=1pt]
     \draw (-0.5,-.0) -- (-0.5,-.4);
     \draw (-10mm,-.9) arc (180:0:5mm);
     \draw (-13mm,-11.75mm) arc (180:0:2.7mm);
     \draw (-3mm,-11.75mm) arc (180:0:2.7mm);
     \draw (-.76,-1.2) -- (-.76,-1.7);
     \draw (-.30,-1.2) -- (-.30,-1.7);
     \draw (-7.7mm,-16.2mm) arc (180:360:2.4mm);
     \draw (-1.32,-1.2) -- (-1.32,-2.1);
     \draw (.23,-1.2) -- (.23,-2.1);
     \filldraw (-.5,-1.85) circle (3pt);
\draw (2.2,-.8) node(x)  {or};
 \draw (4.5,-0) -- (4.5,-.26);
     \draw (4.,-.7) -- (4.,-1.3);
     \draw (39.7mm,-7.5mm) arc (180:0:5mm);
     \draw (47mm,-10.8mm) arc (180:0:3mm);
     \draw (44.5mm,-13.8mm) arc (180:0:3mm);
     \draw (39.9mm,-13.2mm) arc (180:360:2.3mm);
     \filldraw (4.2,-1.6) circle (3pt);
     \draw (5.05,-1.28) -- (5.05,-2.1);
     \draw (5.3,-1.05) -- (5.3,-2.1);
\end{tikzpicture}\end{centering}\caption{.}\end{figure}


\begin{lemma} Let $a$ and $b$ be unaries (modal). Then binary co-operation
$(a,b)$ is coassociative iff $a$ and $b$ are commuting idempotents.\end{lemma}

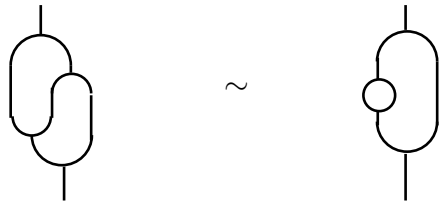
\begin{figure}[!ht]
 \begin{centering}
   \begin{tikzpicture}[line width=1pt]
     \draw (-0.5,-.0) -- (-0.5,-.4);
     \draw (-9mm,-.8) arc (180:0:4mm);
     \draw (-.1,-.8) -- (-.1,-.9);
     \draw (-3.5mm,-11.75mm) arc (180:0:2.6mm);
     \draw (-.9,-.8) -- (-.9,-1.5);
     \draw (-.35,-1.1) -- (-.35,-1.5);
     \draw (-8.75mm,-14.8mm) arc (180:360:2.55mm);
     \draw (.17,-1.2) -- (.17,-1.8);
     \draw (-6.2mm,-17.3mm) arc (180:360:4mm);
     \draw (-.19,-2.15) -- (-.19,-2.6);
     
     \draw (2.1,-1.1) node(x) {$\sim$};
  
     \draw (4.35,-0) -- (4.35,-.33);
     \draw (3.98,-.7) -- (3.98,-.99);
     \draw (39.7mm,-7.5mm) arc (180:0:4mm);
     \draw (4.,-1.2) circle (6pt);
     \draw (39.7mm,-15.5mm) arc (180:360:4mm);
     \draw (3.98,-1.4) -- (3.98,-1.6);
     \draw (4.35,-1.98) -- (4.35,-2.6);
     \draw (4.77,-.75) -- (4.77,-1.6);\end{tikzpicture}
  \end{centering}
  \caption{Co-associativity of binary procreation.}
  \label{coasso}
\end{figure}

Composition of the `elementary' plants from Figures \ref{music} and \ref{binope} looks
like the grafting in an orchard and this grafting is generatig the
operads of composed plants. For example, each of the two grafted plants
$\in\cat(4,2)$ and $\in\cat(2,4)$ on Figure \ref{prodbin} is the results of the grafting of three elementary colored plants from Figure \ref{binope}.

\begin{figure}[!ht]
 \begin{centering}
   \begin{tikzpicture}[line width=1pt]
     \draw (.01,-0.) -- (.01,-.8);
     \draw (.78,-0) -- (.78,-.8);
     \draw (-0mm,-7.8mm) arc (180:360:4mm);
     \draw (.4,-1.18) -- (.4,-2.1);
     \filldraw (.4,-1.6) circle (3pt);
    
     \draw (4.01,-0.) -- (4.01,-.8);
     \draw (4.8,-0) -- (4.8,-.8);
     \draw (40mm,-7.8mm) arc (180:360:4mm);
     \draw (4.39,-1.18) -- (4.39,-2.1);
     \filldraw (4,-.45) circle (3pt);

     \draw (8.01,-0.) -- (8.01,-.8);
     \draw (8.8,-0) -- (8.8,-.8);
     \draw (80mm,-7.8mm) arc (180:360:4mm);
     \draw (8.39,-1.18) -- (8.39,-2.1);
     \filldraw (8.8,-.45) circle (3pt);
     \draw (8,-.45) circle (3pt);
\end{tikzpicture}\end{centering}\caption{Examples of $2\mapsto 1$ grafted with $1\mapsto 1.$}\label{exgraf}\end{figure}
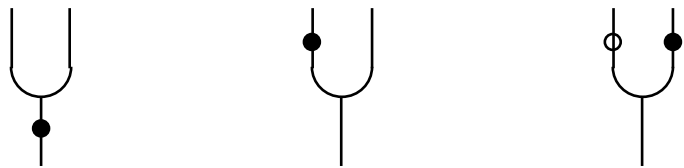

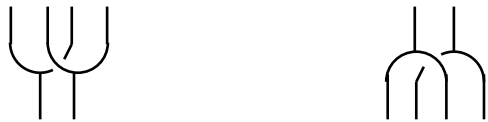
\begin{figure}[!ht]
 \begin{centering}
   \begin{tikzpicture}[line width=1pt]
     \draw (.01,-0.3) -- (.01,-.8);
     \draw (.82,-0.3) -- (.82,-.8);
     \draw (.82,-0.8) -- (.72,-1.);
     \draw (-0mm,-7.8mm) arc (180:295:4mm);
     \draw (.4,-1.18) -- (.4,-1.8);
     \draw (.5,-0.3) -- (.5,-.8);
     \draw (1.29,-0.3) -- (1.29,-.8);
     \draw (5mm,-7.8mm) arc (180:360:4mm);
     \draw (.85,-1.18) -- (.85,-1.8);

     \draw (5.02,-1.2) -- (5.02,-1.8);
     \draw (5.8,-1.2) -- (5.8,-1.8);
     \draw (50mm,-13mm) arc (180:0:4mm);
     \draw (5.38,-.3) -- (5.38,-.9);
     \draw (5.9,-.3) -- (5.9,-.9);
     \draw (63mm,-13mm) arc (0:115:4mm);
     \draw (5.5,-1.1) -- (5.4,-1.3);
     \draw (5.4,-1.3) -- (5.4,-1.8);
     \draw (6.3,-1.2) -- (6.3,-1.8);
\end{tikzpicture}\end{centering}\caption{Product of binary operations and of co-operations, needs the grafting of three plants including braid.}\label{prodbin}\end{figure}

\section{Clone = operad of plants. Hypervariety, Bootstrap}
The grafting of plants is our basic `operation'. Grafting is
multivalued, therefore is not an operation in the usual meaning. From
two copies of binary plant grafting produce two different ternaries,
etc. We prefer use `plant' instead of `operation' because a realization
in sets is not yet assumed (not yet carriers). An algebra in universal
algebra is a set (carrier of algebra) with a family of operations.
Instead, in this paper, we deal with the family of plants of the given
type $<\ldots>$ - a generator of a free operad - and with relations, however
carrier was not yet selected. Therefore by an algebra we mean a family
of plants of the given type $<\ldots>.$

We identify variety of algebras with the quotient operad.

Contrary to distinguished unique co-operation of duplication
$a\mapsto(a,a)$ (Figure \ref{mitosis}) there are no distinguished binary
operation. Figure \ref{mitosis} readed from the bottom define binary with
neutral. Not every binary operation can be constructed in terms of some
distinguished binary and the set of all modal operations.

\begin{exer} For the case of the two element set $\{$true, false$\},$
the classical clone is one generated clone by a $\{$Sheffer stroke$\},$\end{exer}

\begin{figure}[!ht]\begin{centering}
   \begin{tikzpicture}[line width=1pt]
     \draw (-0.5,.1) -- (-0.5,-.3);
     \draw (-.5,-.5) circle (7pt);
     \draw (-0.5,-.8) -- (-0.5,-1.1);
     \draw (-.5,-1.4) circle (7pt);
     \draw (-.5,-1.6) -- (-.5,-2.1);
\draw (1.85,-.9) node{$\sim$};
\draw (3.35,-0) -- (3.35,-2.1);
     \draw [style=dashed](-1.1,-.5) rectangle (.2,-1.4);   
   \end{tikzpicture}
  \end{centering}
  \caption{The Sheffer stroke.}
  \label{sheff}
\end{figure}
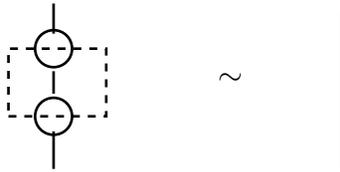

Above we defined an operad on generators (alphabet of morphisms) and relations. This was presentation-dependent quotient operad. One can try to define a clone without of presentation, as presentation-free operad by means of the axiomatic abstract properties of the given operad. We are not going to this subject in the present paper.

\section{The Boolean and Heyting operations}  The \textit{logic} in the present paper is an abbreviation for the algebraic categorical logic. For example the modal logic needs full set of unary operations and their mutual inter-relations.

We define the Boolean operation and the Heyting operation [Heyting 1930, Henkin 1950] as the
relation on two nodes graph with one bubble, as shown on Figure \ref{semilat}.
The Heyting operation is an algebra generated by two plants of type
$\in\cat(2,1)$ and $\in\cat(0,1)$ with this one relation. For example, the Heyting operation enter into BCC-algebra and BCK-algebra invented as the models for implicational
propositional calculus [Is\'eki 1966].

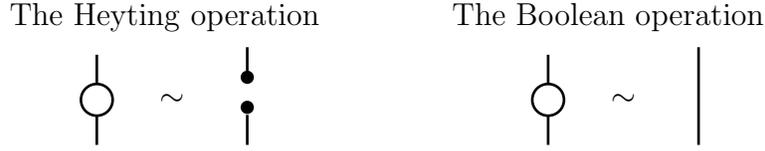
\begin{figure}[!ht]\begin{centering}
   \begin{tikzpicture}[line width=1pt]
     \draw (-0. ,0) -- (-0. ,-.4);
     \draw (-0. ,-.6) circle (6pt);
     \draw (-0. ,-.8) -- (-0. ,-1.2);
     
     \draw (1.0,- .6) node(x)  {$\sim$};

     \draw ( 2. ,.1) -- ( 2. ,-.3);
     \filldraw ( 2. ,-.3) circle (2pt);
     \filldraw ( 2. ,-.7) circle (2pt);
     \draw ( 2. ,-.8) -- ( 2. ,-1.2);
   
     \draw (7.0,- .6) node(x)  {$\sim$};
     \draw ( 6. ,0) -- ( 6. ,-.4);
     \draw ( 6. ,-.6) circle (6pt);
     \draw ( 6. ,-.8) -- ( 6. ,-1.2);
     \draw (.90, .5) node(x)  {The Heyting operation};
     \draw ( 6.8, .5) node(x)  {The Boolean operation};
     \draw (8. ,.1) -- (8. ,-1.2);
   \end{tikzpicture}
  \end{centering}
  \caption{Two semi-lattices: the Heyting semi-lattice and the Boolean semi-lattice.}
  \label{semilat}
\end{figure}

The first node $\in\cat(1,2)$ represents either the pure duplication or a duplication
with essential mutants (an emission, pro-creation, decomposition of
information, etc) and the next node $\in\cat(2,1)$ represents an essential binary
(primitive, \ie not unary with killer), absorption, anihilation,
consumption, reinforcement, $\ldots.$

We do not assume that the Boolean and the Heyting binaries need to be
associative and/or commutative. If instead for some finite integer $n>1$
the (abstract) minimum polynomial of this bubble is $B^n=\id$ (\ie the
bubble $B$ is neither Boolean, nor the Heyting), then is said to be
(generalized) Shefferian, Pierce, or Post operation.

The Heyting operation generalize an implicative algebra [Rasiowa 1974,
p. 16]. An implicative algebra play for some non-classical logics a role
analogous to that played by the Boolean algebra for the classical logic.

For the Boolean semilattice, the Heyting semilattice, the Boolean
$\Z_2$-algebra, the Heyting algebra, we refer to [Lambek \@ Scott, p.
36, Examples 7.3-7.4].

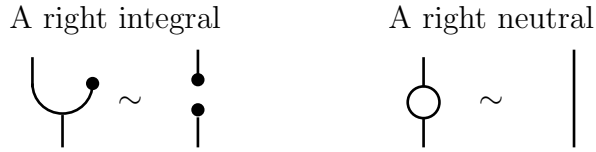
\begin{figure}[!ht]
 \begin{centering}
   \begin{tikzpicture}[line width=1pt]
     \draw (-.2 ,0) -- (-.2 ,-.4);
     \filldraw (.6 ,-.35) circle (2pt);
     \draw (-2mm,-3.5mm) arc (180:360:4mm);
     \draw (0.2 ,-.77) -- (0.2 ,-1.2);
     
     \draw (1.1,- .6) node(x) {$\sim$};

     \draw ( 2. ,.1) -- ( 2. ,-.3);
     \filldraw ( 2. ,-.3) circle (2pt);
     \filldraw ( 2. ,-.7) circle (2pt);
     \draw ( 2. ,-.8) -- ( 2. ,-1.2);
   
     \draw (5.9,- .6) node(x) {$\sim$};
     \draw ( 5. ,0) -- ( 5. ,-.4);
     \draw ( 5. ,-.6) circle (6pt);
     \draw ( 5. ,-.8) -- ( 5. ,-1.2);
     \draw ( .9, .50) node(x)  {A right integral};
     \draw ( 5.9, .50) node(x)  {A right neutral};
     \draw (7. ,.1) -- (7. ,-1.2);
   \end{tikzpicture}
  \end{centering}
  \caption{Integral and neutral.}
  \label{intneut}
\end{figure}

\section{Lattice and Quasi-Group}The filled and not filled circles on Figure \ref{quasilat} represent two binaries. At this general setting neither associativity nor commutativity need to
be assumed. We are contrasting the interrelations between two binaries defining the (right) lattice (where the mather survive) and the (right) quasigroup (where mother is cancelled [Moufang 1935]). Both relations (equations, identities) on Figure \ref{quasilat} define \textit{strongly nonregular} clones (or varieties) considered by Graczy\'nska [1989, 1990, 1998, 2006]. Terminology \textit{regular} clone, etc, was introduced by [P{\l}onka 1969].

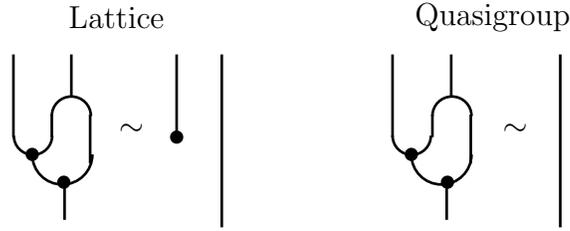
\begin{figure}[!ht]
 \begin{centering}
   \begin{tikzpicture}[line width=1pt]
     \draw (- .6,  .3) -- (- .6, -.25);
     \draw (-8.6mm,-5.2  mm) arc (180:0:2.6mm);
     \draw (-1.37, .3) -- (-1.37,- .8);
     \draw (-.86,- .5) -- (-.86,- .8);
     \draw (-13.65mm,-7.8mm) arc (180:360:2.55mm);
     \filldraw (-1.12 ,-1.03) circle (2pt);
     \draw ( -.35 ,- .5) -- (-.35,-1.15);
     \draw (-11.2mm,-10.3mm) arc (180:360:4mm);
     \filldraw (-.7 , -1.4) circle (2pt);
     \draw (-.69, -1.4) -- (-.69,-1.9);
     
     \draw ( .2,-.7) node(x)  {$\sim$};
     
     \draw (.8 ,.3) -- (.8 ,- .8);
     \draw (1.4 ,.3) -- (1.4,-2);
     \filldraw ( .8, -.8) circle (2pt);
     
     \draw ( 4.45,  .3) -- ( 4.45, -.25);
     \draw (  42  mm,-5.2  mm) arc (180:0:2.6mm);
     \draw ( 3.67, .3) -- ( 3.67,- .8);
     \draw ( 4.20,- .5) -- ( 4.20,- .8);
     \draw ( 36.7mm,-7.8mm) arc (180:360:2.55mm);
     \filldraw ( 3.92 ,-1.03) circle (2pt);
     \draw (4.71,- .5) -- (4.71,-1.15);
     \draw ( 39.2mm,-10.3mm) arc (180:360:4mm);
     \filldraw ( 4.4 , -1.4) circle (2pt);
     \draw ( 4.39, -1.4) -- ( 4.39,-1.9);
     
     \draw (  5.3,-.7) node(x)  {$\sim$};
     \draw (  .0,  .8) node(x)  {Lattice};
     \draw ( 5. 0,  .8) node(x)  {Quasigroup};
     
     \draw (6.5 ,.3) -- (6.5 ,-1.1);
     \draw (5.9 ,.3) -- (5.9 ,-2);
     \filldraw ( 6.5, -1.1) circle (2pt);
\end{tikzpicture}\end{centering}
  \caption{The right lattice and the right quasigroup (the right cancellation).}
  \label{quasilat}
\end{figure}

\section{Braid} An endomap of arity (type) $\in\cat(2,2),$ is said
to be \textit{pre-braid} if the Artin prebraid relation $\in\cat(3,3),$
represented by tangles on Figure \ref{artbraid}, holds. An invertible prebraid is
said to be a braid.
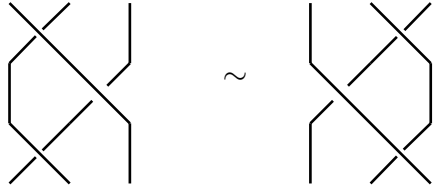
\begin{figure}[!ht]
 \begin{centering}
   \begin{tikzpicture}[line width=1pt]
     \draw (0.,  0) -- (1.6,-1.6);
     \draw (.8 , .0) -- (0.44,-.35);
     \draw (0., -.8) -- (.35,-.44);
     \draw (1.6 , .0) -- (1.6,-.8);
     
     \draw (0,-.8) -- (0,-1.6);
     \draw (1.6,- .8) -- (1.3,-1.1);
     \draw (.8,-1.6) -- (1.1,-1.3);
     
     \draw (0.,-1.6) -- (.8,-2.4);
     \draw (.8 ,-1.6) -- (0.44,-1.96);
     \draw (0.,-2.4) -- (.35,-2.05);
     \draw (1.6 ,-1.6) -- (1.6,-2.4);
    
     \draw (  3. ,-1.0) node(x)  {$\sim$};
 \draw (4.  , .0) -- (4. ,-.8);
     \draw (4.8,-0) -- (5.6,- .8);
     \draw (5.6 ,-0) -- (5.25,-.36);
     \draw (4.8,-.8) -- (5.15,-.45);
     
     \draw (4.8,- .8) -- (4.5,-1.1);
     \draw (4. ,-1.6) -- (4.3,-1.3);
     \draw (5.6 ,- .8) -- (5.6,-1.6);
     
     \draw (4.  ,-1.6) -- (4. ,-2.4);
     \draw (4. ,-.8) -- (5.6,-2.4);
     \draw (5.6 ,-1.6) -- (5.24,-1.95);
     \draw (4.8, -2.4) -- (5.15,-2.04);
 \end{tikzpicture}\end{centering}
  \caption{The Artin prebraid relation.}\label{artbraid}\end{figure}

Let plant $\in\cat(2,2)$ is composed from a pair of from $\alpha,\beta\in\cat(2,1).$
Therefore we can insert into every vertex on Figure \ref{artbraid} the brassiere Figure \ref{brass}, \ie a pair $(\alpha,\beta)$ of the binary operations
\begin{figure}[!ht]
 \begin{centering}
   \begin{tikzpicture}[line width=1pt]
     \draw (-2.5,-.4) -- (-1.,-1.9);
     \draw (-1. ,-.4) -- (-1.6,-1. );
     \draw (-2.5,-1.9) -- (-1.9,-1.3);

     \draw ( 1.0,-1.3) node(x)  {$\leadsto$};

     \draw (3.2,-.4) -- (3.2  ,- .7);
     \draw (3.2,-.7) -- (4.,-1.6);
     \draw (4. ,-.7) -- (3.7,-1.1);
     \draw (3.2,-1.6) -- (3.5,-1.25);
     \draw (32mm,- 7.3mm) arc ( 90:270:4.4mm);
     \draw (40mm,- 7.3mm) arc (270:90:-4.4mm);
     \draw (4.0,-.4) -- (4.0,- .7);
     \draw (  3.0,-1.8) node(x)  {$\alpha$};
     \draw (  4.25,-1.8) node(x)  {$\beta$};
     \draw (3.2,-1.6) -- (3.2,-1.9);
     \draw (4. ,-1.6) -- (4. ,-1.9);
\end{tikzpicture}
  \end{centering}
  \caption{Brassiere.}
  \label{brass}
 \end{figure}
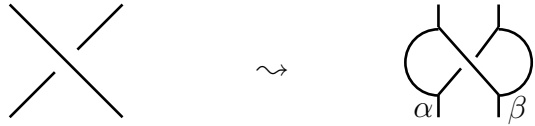

Then it is easily seen that the Artin relation on $3\mapsto 3$ is
equivalent to three ternary (regular? [P{\l}onka]) relations shown on
Figure \ref{artpreb}.

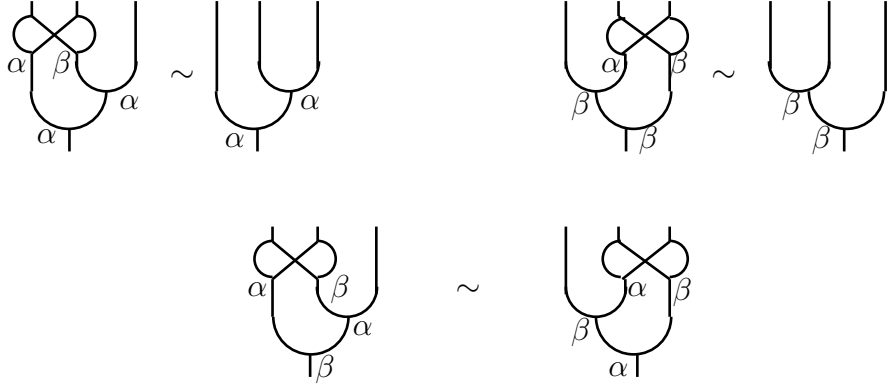
\begin{figure}[!ht]
 \begin{centering}
   \begin{tikzpicture}[line width=1pt]
     \draw (-1.3,-.4) -- (-1.3 ,- .6);%
     \draw (-1.3,-.6) -- (-1.9,-1.1);
     \draw (-1.9 ,-.6) -- (-1.3,-1.1);
     \draw (-13mm,- 10.8mm) arc ( 90:270:-2.4mm);
     \draw (-19mm,- 10.8mm) arc (270:90:2.4mm);
     \draw (-1.9,-.4) -- (-1.9,- .6);
     \draw (- 2.1,-1.25) node(x)  {$\alpha$};
     \draw (- 1.5 ,-1.25) node(x)  {$\beta$};
     \draw (-  .6,- 1.75) node(x)  {$\alpha$};
     \draw (-1.7,-2.2) node(x)  {$\alpha$};
     \draw (  .1 ,-1.4) node(x)  {$\sim$};
     \draw (  1.8 ,-1.75) node(x)  {$\alpha$};
     \draw (  .8 ,-2.25) node(x)  {$\alpha$};

     \draw (-1.9,-1.1) -- (-1.9,-1.6);
     \draw (-1.3 ,-1.1) -- (-1.3 ,-1.3);
     
     \draw (- .52,- .4) -- (- .52,-1.3);
     \draw (-13.1mm,-12. mm) arc ( 180:360:4 mm);
     \draw (-19.1mm,-16. mm) arc ( 180:360:5 mm);
     \draw (-1.4 ,-2.1) -- (-1.4 ,-2.4);

     \draw ( 1.13,-.4) -- ( 1.13,-1.3);%
     \draw (1.9,-.4) -- (1.9,- 1.3);
     \draw (  .56 ,- .4) -- (  .56 ,-1.6);
     \draw ( 11.18mm,-12. mm) arc ( 180:360:4 mm);
     \draw (  5.5mm,-16. mm) arc ( 180:360:5 mm);
     \draw (1.1 ,-2.1) -- (1.1 ,-2.4);
     \draw (5.9,-.4) -- (5.9  ,- .6);%
     \draw (5.9,-.6) -- (6.58,-1.1);
     \draw (6.58 ,-.6) -- (5.95,-1.1);
     \draw (59.9mm,- 6.3mm) arc ( 90:270:2.4mm);
     \draw (65.8mm,- 6.3mm) arc (270:90:-2.4mm);
     \draw (6.58,-.4) -- (6.58,- .6);
     \draw (5.99,-1.1) -- (5.99,-1.3);
     \draw (6.58 ,-1.1) -- (6.58 ,-1.6);
     
     \draw (5.2  ,- .4) -- (5.2  ,-1.3);
     \draw (52  mm,-12. mm) arc ( 180:360:4 mm);
     \draw (56 mm,-16. mm) arc ( 180:360:5 mm);
     \draw (6.0 ,-2.1) -- (6.0 ,-2.4);
     \draw ( 5.8,-1.25) node(x)  {$\alpha$};
     \draw ( 6.7,-1.25) node(x)  {$\beta$};
     \draw ( 5.4,-1.8 ) node(x)  {$\beta $};
     \draw ( 6.3,-2.25) node(x)   {$\beta $};
     \draw ( 7.3,-1.4) node(x)  {$\sim$};
     \draw ( 8.2 ,-1.8 ) node(x) {$\beta$};
     \draw ( 8.6,-2.25) node(x)  {$\beta$};

     \draw (8.7,-.4) -- (8.7  ,-1.3);%
     \draw (9.45 ,- .4) -- (9.45 ,-1.6);
     \draw (7.92,-.4) -- (7.92  ,-1.3);%
     \draw (79mm,-12. mm) arc ( 180:360:4 mm);
     \draw (84.3mm,-16. mm) arc ( 180:360:5 mm);
     \draw (8.9 ,-2.1) -- (8.9 ,-2.4);

     \draw ( 1.3,-3.4) -- ( 1.3 ,-3.6);%
     \draw ( 1.3,-3.6) -- ( 1.9,-4.1);
     \draw ( 1.9 ,-3.6) -- ( 1.3,-4.1);
     \draw ( 13mm,- 35.9mm) arc ( 90:270: 2.4mm);
     \draw ( 19mm,- 35.9mm) arc (270:90:- 2.4mm);
     \draw ( 1.9,-3.4) -- ( 1.9,- 3.6);
     \draw ( 1.3,-4.1) -- ( 1.3,-4.6);
     \draw ( 1.9 ,-4.1) -- ( 1.9 ,-4.3);
     
     \draw ( 2.68,-3.4) -- ( 2.68,-4.3);
     \draw ( 18.9mm,-42. mm) arc ( 180:360:4 mm);
     \draw ( 13.1mm,-46. mm) arc ( 180:360:5 mm);
     \draw ( 1.8 ,-5.1) -- ( 1.8 ,-5.4);

     \draw (5.9,-3.4) -- (5.9  ,-3.6);%
     \draw (5.9,-3.6) -- (6.58,-4.1);
     \draw (6.58 ,-3.6) -- (5.95,-4.1);
     \draw (59.3mm,-35.9mm) arc ( 90:270:2.4mm);
     \draw (65.8mm,-35.9mm) arc (270:90:-2.4mm);
     \draw (6.58,-3.4) -- (6.58,- 3.6);
     \draw (5.99,-4.1) -- (5.99,-4.3);
     \draw (6.58 ,-4.1) -- (6.58 ,-4.6);
     \draw (  1.1,-4.23) node(x)  {$\alpha$};
     \draw (  2.2,-4.23) node(x)  {$\beta$};
     \draw (  2.5 ,-4.75 ) node(x)  {$\alpha$};
     \draw (  2.0,-5.3 ) node(x)  {$\beta $};
     \draw (  3.9 ,-4.25) node(x)  {$\sim$};
     \draw (  6.15,-4.23) node(x)  {$\alpha$};
     \draw (  6.75,-4.27) node(x)  {$\beta$};
     \draw (   5.4,-4.8 ) node(x)  {$\beta  $};
     \draw (  5.9,-5.3) node(x)  {$\alpha$};

     \draw (5.2  ,-3.4) -- (5.2  ,-4.3);
     \draw (52  mm,-42. mm) arc ( 180:360:4 mm);
     \draw (56 mm,-46. mm) arc ( 180:360:5 mm);
     \draw (6.15,-5.1) -- (6.15,-5.4);
\end{tikzpicture}
  \end{centering}
  \caption{A pair of binary operation as the Artin prebraid.}
  \label{artpreb} \end{figure}

\begin{defC}[Essential morphism] All operations $\in\cat(m,1)$ (and also all plants $\in\cat(m,n)$) are segregated into two groups: not essential (or trivial) and essential or
primitives (not trivials). An unary operation $\in\cat(1,1)$ is said to
trivial if it is a constant map (the composition of the killer with
nullary operation) or the identity plant. All other unaries are said to
be essentials.\end{defC} 

For example in realization on the finite set $S$ with
the cardinality $s\equiv|S|,$ the set of all unaries has the
cardinality $s^s.$ The number of primitive unaries is
$$s^s-s-1=\begin{cases}1&\text{if $s=2$ \ie in classical logic,}\\
23&\text{if $s=3$ in the {\L}ukasiewicz logic.}\end{cases}$$

A plant $1<m\mapsto 1$ is said to be derived (or grafted, not essential,
trivial) if can be build by juxtapositioning and grafting from plants
$m>k\mapsto 1.$ For example the number of the essential binaries
$2\mapsto 1$ in the realization on the finite set $S$ as above, is
$$s^{(s^2)}-2s^s+s=\begin{cases}2^4-6=10&\text{in the classical case,}
\\3^9-51=19632&\text{in the {\L}ukasiewicz logic.}\end{cases}$$

Therefore among 16 binary classical connectives, 6 are not primitives
because they are juxtaposition of unary with killer, and 10 are
essential. In particular `projector' $m\mapsto 1$ is juxtaposition of
the killers with the identity plant.

Our first result concern the classical two-valued logic. In this case
a priori there are 100 pairs of primitives binaries as the candidates
for the Artin prebraid.

\begin{cor} Let in Figure \ref{artpreb} the both binary plants be
primitive. Then in the classical (two-valued logic) there are four
prebraids build from the disjunction and the conjunction only. All these
prebraids are idempotents.\end{cor}

\begin{cor} In the classical logic both the disjunction and the
conjunction are $(\alpha,\beta)$-symmetric.\end{cor}
\begin{proof}This statement is equivalent to the idempotency of the pre-braids.\end{proof}

If in a plant $\in\cat(2,2),$ $\alpha\equiv\beta,$ then brassiere is
reduced to the duplication of one binary as shown on the first graph in
Figure \ref{artpreb}. Then every boolean operation (necessarily associative)
gives prebraid.

The other four graphs on Figure \ref{artpreb} correspond to the cases when at
least one among two binaries in the brassiere is given by unary (with
killer). In these cases a binary cooperation needs not to be just
duplication.
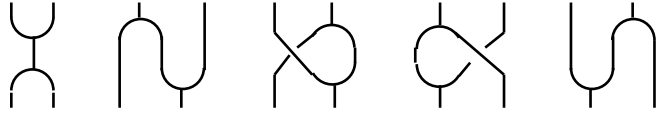
\begin{figure}[!ht]
 \begin{centering}
   \begin{tikzpicture}[line width=1pt]
     \draw ( 1.9,-.4) -- ( 1.9 ,-.6);%
     \draw ( 2.46,-.4) -- ( 2.46 ,-.6);%
     \draw ( 19mm,- 5.9 mm) arc (180:360: 2.8mm);
     \draw ( 2.2,- .85) -- ( 2.2,- 1.3);
     \draw ( 24.6mm,- 15.7mm) arc ( 0:180: 2.8mm);
     \draw ( 1.9,-1.6) -- ( 1.9 ,-1.8);%
     \draw ( 2.46,-1.6) -- ( 2.46 ,-1.8);%
     
     \draw ( 3.6,-.4) -- ( 3.6 ,-.6);%
     \draw ( 4.47,-.4) -- ( 4.47  ,-1.3);%
     \draw ( 39.mm,- 12.7 mm) arc (180:360: 2.8mm);
     \draw ( 3.9,- .85) -- (3.9,- 1.3);
     \draw ( 39mm,- 9.mm) arc ( 0:180: 2.8mm);
     \draw ( 3.34,-.85) -- ( 3.34 ,-1.8);%
     \draw ( 4.16,-1.56) -- ( 4.16 ,-1.8);%
     
     \draw ( 6.16,-.4) -- ( 6.16 ,-.7);%
     \draw ( 6.47,-1.) -- ( 6.47,-1.2);%
     \draw ( 59.1mm,-13.4 mm) arc (210:360: 3. mm);
     \draw ( 5.4,- .4) -- ( 5.4 ,- .8);%
     \draw ( 5.4,- .8) -- ( 5.9,- 1.36);
     \draw ( 5.95,- .8) -- ( 5.7,- 1.0);
     \draw ( 5.4,- 1.36) -- ( 5.6,- 1.1);
     \draw ( 5.4,-1.36) -- ( 5.4 ,-1.8);%
     \draw ( 64.6mm,- 9.99mm) arc ( 0:150: 3. mm);
     \draw ( 6.2,-1.5) -- ( 6.2 ,-1.8);%
     
     \draw ( 7.6 ,-.4) -- ( 7.6  ,-.7);%
     \draw ( 7.27,-1.) -- ( 7.27,-1.2);%
     \draw ( 72.9mm,-12.2 mm) arc (180:330: 3. mm);
     \draw ( 8.45,- .4) -- ( 8.45 ,- .8);%
     \draw ( 7.8,- .8) -- ( 8.45,- 1.36);
     \draw ( 8.45,- .8) -- ( 8.2,- 1.0);
     \draw ( 7.85,- 1.38) -- ( 8.0,- 1.2);
     \draw ( 8.45,-1.36) -- ( 8.45 ,-1.8);%
     \draw ( 78.5mm,- 8.69mm) arc (30:180: 3. mm);
     \draw ( 7.6,-1.5) -- ( 7.6 ,-1.8);%
     
     \draw ( 9.6,-1.56) -- ( 9.6 ,-1.8);%
     \draw (10.45,-.9) -- (10.45  ,-1.8);%
     \draw ( 93.4mm,- 12.7 mm) arc (180:360: 2.8mm);
     \draw ( 9.9,- .85) -- (9.9,- 1.3);
     \draw (104.5mm,- 9.mm) arc ( 0:180: 2.8mm);
     \draw ( 9.34,-.4) -- ( 9.34 ,-1.3);%
     \draw (10.16,- .4) -- (10.16 ,-.6);\end{tikzpicture}\end{centering}
  \caption{Examples of grafted $\in\cat(2,2).$}
 \label{exgraf2}\end{figure}

\begin{conjecture} If at least one binary among $(\alpha,\beta)$ in Figure \ref{artpreb}
 is unary (with killer) then the last Artin relation in Figure \ref{artpreb}
is equivalent that $\alpha\in\hom\beta$ or $\beta\in\hom\alpha.$
\end{conjecture}

Figure \ref{binmor} gives the two examples of the `identities', two examples of
the possible laws of the nature. Two plants grafted on the left of each
law are identical as three plants grafted on the right. These laws are
interpreted that a binary $\in\cat(2,1)$ (an operation, an action, et
cetera) is the morphism with respect to the process $\in\cat(2,2),$
the left law is \textit{under} morphism, the right law is \textit{over}
morphism. In particular Figure \ref{artbraid} is nothing but morphism
(simultaneously under \textit{and} over) of a prebraid $\in\cat(2,2)$ with
respect to himself, \ie selfmorphism.
\begin{figure}[!ht]
 \begin{centering}
   \begin{tikzpicture}[line width=1pt]
     \draw (- 4.,   .5) -- (- 4.0, - .3);
     \draw (-3.29,   .5) -- (- 3.29,   .3);
     \draw (-2.69,  .5) -- (-2.69,   .3);
     \draw (-33. mm,  3.4  mm) arc (180:360:3mm);
     \draw (- 3.  ,  .05) -- (- 3.  ,-  .3);
     \draw (-2.0,-.5) node(x)  {$\sim$};
     \draw ( 4.8,-.5) node(x)  {$\sim$};
     \draw (-4.0,-  .3) -- (-3.  ,- 1.4);
     \draw (- 3. ,-  .3) -- (-3.4,-  .8);
     \draw (-4.0,- 1.4) -- (-3.6,- 1. );
     
     \draw (- 1.0,  .5) -- ( .8,- .8);
     \draw (- .10,  .5) -- (- .5 , .2);
     \draw (-1 .0,- .1) -- ( - .70,  .1 );
     \draw ( .80,   .5) -- ( -0.00, - .1);
     \draw (- .42,- .4) -- ( -.21,-  .25);
     \draw (-1.03  ,-  .1) -- ( -1.03  ,-  .6);
     \draw (- .42,- .4) -- (- .42,-  .7);
     \draw (-10.3 mm,- 5.7  mm) arc (180:360:3.05mm);
     \draw (- .72,- .9) -- (- .72,-1.4);
     \draw ( .80,- .8) -- ( .80,-1.4);
     \draw (  4.,   .5) -- (  4.0, - .3);
     \draw ( 3.29,   .5) -- (  3.29,   .3);
     \draw ( 2.69,  .5) -- ( 2.69,   .3);
     \draw ( 26.9mm,  3.4  mm) arc (180:360:3mm);
     \draw (  3.  ,  .05) -- (  3.  ,-  .3);
     \draw ( 3.0,-  .3) -- ( 4.  ,- 1.4);
     \draw (  4. ,-  .3) -- ( 3.6,-  .8);
     \draw ( 3.0,- 1.4) -- ( 3.4,- 1. );
     
     \draw (5.80,   .5) -- (  6.80, - .35);
     \draw (  6.4,  .5) -- (7.41,- .3);
     \draw ( 7.41,  .5) -- ( 6.9 , .2);
     \draw ( 6 .5, .0) -- (  6.70,  .1 );
     \draw ( 5.80,- .4) -- ( 6.3,-  .1);
      \draw ( 6.80  ,-  .35) -- (  6.80  ,-  .6);
      \draw ( 7.41,- .3 ) -- ( 7.41,-  .6);
     \draw ( 68. mm,- 5.7  mm) arc (180:360:3.05mm);
     \draw ( 7.12,- .9) -- ( 7.12,-1.4);
     \draw (5.80,- .4) -- (5.80,-1.4);
   \end{tikzpicture}
  \end{centering}
  \caption{A binary as the morphism: under and over.}
  \label{binmor}
\end{figure}
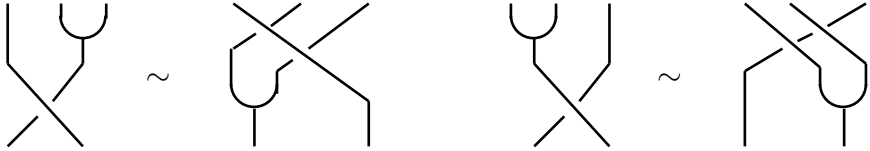

Specializing Figure \ref{binmor} to the relations expressing that binary
$\gamma$ is the $(\alpha,\beta)$-morphism we get the system of
(regular?) relations on Figures \ref{gamunmor} and \ref{gamovmor}.

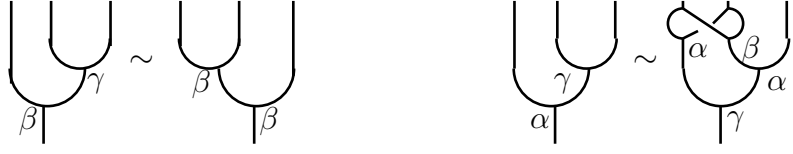
\begin{figure}[!ht]
 \begin{centering}
   \begin{tikzpicture}[line width=1pt]
     \draw (- 1.83, 1.3) -- (- 1.83,  .20);%
     \draw (- 1.3, 1.3) -- (-1.3,   .7);
     \draw (  -.51 , 1.3) -- (- .51 ,  .7);
     \draw (-  1.6,- .3 ) node(x)  {$\beta$};
     \draw (-  .7, .2 ) node(x)  {$\gamma $};
     \draw (-  .1 ,  .5) node(x)  {$\sim$};
     \draw (   .7,  .2) node(x)  {$\beta$};
     \draw (  1.6, - .3) node(x)  {$\beta$};
     \draw (-13.18mm,  8. mm) arc ( 180:360:4 mm);
     \draw (-18.5mm,  4.  mm) arc ( 180:360:5 mm);
     \draw (- 1.4 , - .1) -- (- 1.4 ,- .6);
     \draw ( 1.2,1.3) -- (1.2,.7);%
     \draw (1.92 ,1.3) -- (1.92 ,.3);
     \draw ( .42,1.3) -- ( .42  ,.7);%
     \draw ( 4mm, 8. mm) arc ( 180:360:4 mm);
     \draw ( 9.3mm,  4. mm) arc ( 180:360:5 mm);
     \draw (1.4 ,- .1) -- (1.4 ,- .6);
      \draw ( 5.43,1.3) -- ( 5.43,  .8);%
     \draw (6.2,1.3) -- (6.2,   .8);
     \draw ( 4.86 , 1.3) -- ( 4.86 ,  .4);
     \draw ( 54.18mm,  8. mm) arc ( 180:360:4 mm);
     \draw ( 48.5mm,  4. mm) arc ( 180:360:5 mm);
     \draw (5.4 ,- .1) -- (5.4 ,- .6);

     \draw ( 7.1,1.2) -- ( 7.7 ,  .8);%
     \draw ( 7.7,1.2) -- ( 7.5, 1. );
     \draw ( 7.3 , .9) -- ( 7.1,  .8);
     \draw ( 77mm,   8. mm) arc ( 90:270:-2. mm);
     \draw ( 71mm,   8. mm) arc (270:90:2. mm);
     \draw ( 7.1,1.3) -- ( 7.1,  1.2);
     \draw ( 7.7, 1.3) -- ( 7.7, 1.2);
     \draw ( 7.1 ,  .8) -- ( 7.1 ,  .4);
     \draw (  5.2,- .3 ) node(x)  {$\alpha$};
     \draw (  5.5,  .2) node(x)  {$\gamma$};
     \draw (  6.6 ,  .5) node(x)  {$\sim$};
     \draw (  7.3,  .65) node(x)  {$\alpha$};
     \draw (  8.0, .65) node(x)  {$\beta$};
     \draw (  8.35 ,  .2) node(x)  {$\alpha$};
     \draw (  7.8,-  .3) node(x)  {$\gamma$};
     
     \draw ( 8.48, 1.3) -- ( 8.48,  .8);
     \draw ( 77.1mm,  8. mm) arc ( 180:360:4 mm);
     \draw ( 71.1mm,  4. mm) arc ( 180:360:5 mm);
     \draw ( 7.6 ,- .1) -- ( 7.6 ,- .6);
  \end{tikzpicture}\end{centering}
  \caption{Binary $\gamma$ is under $(\alpha,\beta)$-morphism.}
  \label{gamunmor} \end{figure}

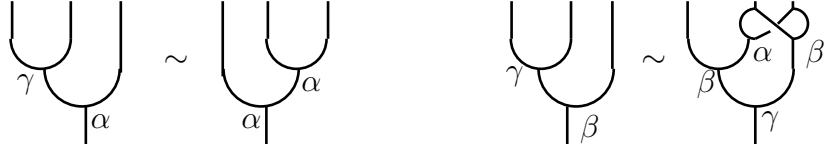
\begin{figure}[!ht]
 \begin{centering}
   \begin{tikzpicture}[line width=1pt]
     \draw (- .535, 1.3) -- (- .535,  .35);%
     \draw (- 1.2, 1.3) -- (-1.2,   .8);
     \draw (  -1.98 , 1.3) -- (-1.98 ,  .8);
     \draw ( - .8,- .3 ) node(x)  {$\alpha$};
     \draw (- 1.8,  .2) node(x)  {$\gamma $};
     \draw (   .2 ,  .5) node(x)  {$\sim$};
     \draw ( 2.0,  .2 ) node(x)  {$\alpha$};
     \draw ( 1.2,-  .3) node(x)  {$\alpha$};
     \draw (-20.mm, 8. mm) arc ( 180:360:4 mm);
     \draw (-15.5mm,  4. mm) arc ( 180:360:5 mm);
     \draw (- 1.  , - .1) -- (- 1.  ,- .6);
     
     \draw (1.42,1.3) -- (1.42,.8);%
     \draw (.82  ,1.3) -- (.82 ,.3);
     \draw (2.22,1.3) -- ( 2.22  ,.8);%
     \draw (  14.  mm,  8. mm) arc ( 180:360:4 mm);
     \draw (  8.3mm,  4.  mm) arc ( 180:360:5 mm);
     \draw (1.4 ,- .1) -- (1.4 ,- .6);  
     \draw ( 4.66,1.3) -- ( 4.66,  .8);%
     \draw (5.45,1.3) -- (5.45,   .8);
     \draw ( 6. , 1.3) -- ( 6. ,  .4);
     \draw ( 46.5mm,  8. mm) arc ( 180:360:4 mm);
     \draw ( 50.18mm, 4. mm) arc ( 180:360:5 mm);
     \draw (5.4 ,- .1) -- (5.4 ,- .6);

     \draw ( 7.9,1.2) -- ( 8.4 ,  .8);%
     \draw ( 8.4,1.2) -- ( 8.2, 1. );
     \draw ( 8.1 , .9) -- ( 7.9,  .8);
     \draw ( 84mm,   8. mm) arc ( 90:270:-2. mm);
     \draw ( 79mm,   8. mm) arc (270:90:2. mm);
     \draw ( 7.9,1.3) -- ( 7.9,  1.2);
     \draw ( 8.4, 1.3) -- ( 8.4, 1.2);
     \draw ( 8.4 ,  .8) -- ( 8.4 ,  .4);
     \draw (  5.7,- .4) node(x)  {$\beta $};
     \draw (  4.7, .3) node(x)  {$\gamma$};
     \draw (  6.55,   .5) node(x)  {$\sim$};
     \draw (  8.  ,  .6 ) node(x)  {$\alpha$};
     \draw (  8.7,  .6) node(x)  {$\beta$};
     \draw (  7.25,  .2 ) node(x)  {$\beta $};
     \draw (  8.1,- .3 ) node(x)  {$\gamma$};
     
     \draw ( 7.  , 1.3) -- ( 7.  ,  .8);
     \draw ( 70.1mm,  8. mm) arc ( 180:360:4 mm);
     \draw ( 74.1mm,  4. mm) arc ( 180:360:5 mm);
     \draw ( 7.9 ,- .1) -- ( 7.9 ,- .6);
\end{tikzpicture}
  \end{centering}
  \caption{Binary $\gamma$ is over $(\alpha,\beta)$-morphism.}
  \label{gamovmor}
 \end{figure}

The relations on the left in Figures \ref{gamunmor} and \ref{gamovmor} coincide with the so called condition $(S)$ in BCK-algebra of a type
\begin{gather}<2,2,0>\simeq\cat(2,1),\cat(2,1),\cat(0,1).\end{gather}

\begin{thm} Let $\alpha$ and $\beta$ be two Heyting binaries and
$(\alpha,\beta)$-mor\-ph\-i\-sms. Then a clon (of abstract algebras)
$(\alpha,\beta,0,1)$ of type $(2,2,0,0)$ is quasigroup.
\end{thm}
\begin{proof} In Figures \ref{gamunmor} and \ref{gamovmor} we must take $\gamma=\alpha$ and
$\gamma=\beta$ and graft the duplication on the tops. Then insert the
Heyting operation from Figure \ref{sheff} with two nullary operations.
\end{proof}

If a plant of cooperation is under $\sigma$-morphism then
the tangles in Figure \ref{binmor} can be reexpressed as braided plant
as shown on Figure  \ref{undermor}.

\begin{figure}[!ht]
 \begin{centering}
   \begin{tikzpicture}[line width=1pt]
     
     \draw ( 1.06, .4) -- ( 1.06 , .15);%
     \draw ( 14.6mm,- 1.99mm) arc ( 5:150: 4. mm);
      \draw (  2.5,- .8) node(x)  {$\sim$};
     \draw ( 1.47,- .2) -- ( 1.47,- 1.2);%
     \draw (   .7,- .05) -- (   .7 ,- .2);%
     \draw (-  .2,  .4) -- (-  .2 ,- .2);%
     \draw (-  .2,- .2) -- (   .7,- 1.2);
     \draw (   .7 ,- .2) -- (   .4,-  .6);
     \draw (-  .2,- 1.2 ) -- (   .2,-  .8);
     \draw (-  .2,-1.2 ) -- (-  .2 ,-1.8);%
     \draw (   .7,-1.2) -- (   .7 ,-1.35);%
     \draw (  7.3mm,-13.4 mm) arc (210:355: 4. mm);
     \draw ( 1.2,-1.5) -- ( 1.2 ,-1.8);%
     
     \draw ( 3.6, .4) -- ( 4.47 ,-.2);%
     \draw ( 3.6,-.2 ) -- ( 3.9 , .0);%
     \draw ( 4.47, 0.4) -- (4.15 , .1);%
     
     \draw ( 3.6,-.2) -- ( 3.6 ,-.4);%
     \draw ( 3.34,-.7 ) -- ( 3.34 ,-1.8);%
     \draw ( 39mm,- 7.mm) arc ( 0:180: 2.8mm);
     \draw ( 4.47,-.2) -- ( 4.47  ,- .7);%
     
     \draw ( 4.47,-.4) -- ( 4.47  ,- .7);%
     \draw ( 3.9,- .7 ) -- (4.47,-1.3 );
     \draw ( 3.9,-1.3 ) -- (4.17,-1.05);
     \draw ( 4.27,- .9 ) -- (4.47,- .7);
     \draw ( 39.mm,- 12.7 mm) arc (180:360: 2.8mm);
     \draw ( 4.16,-1.56) -- ( 4.16 ,-1.8);\end{tikzpicture}
  \end{centering}\caption{Under $\sigma$-morphism.}
 \label{undermor}\end{figure}
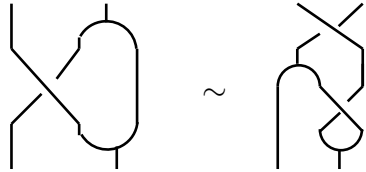

\section{Cooperation as the morphism in a braided category}
When we wish to include an cooperation from Figure  \ref{cooper} as the morphism
in the prebraided monoidal category generated by prebraid
$(\alpha,\beta),$ then this leads naturally to new two candidates for
the Artin prebraids, so to say we get, `the system' of prebraids. Two
new processes $\in\cat(2,2)$ which enter into the game are given on Figure \ref{growsys}.

\begin{figure}[!ht]
 \begin{centering}
   \begin{tikzpicture}[line width=1pt]
\draw (-2.6,1. ) -- (-2.6 , .6);%
     \draw (-1.74,1. ) -- (- 1.74  ,- .1);%
     \draw (-23.mm,-   .7 mm) arc (180:360: 2.8mm);
     \draw (-2.3 ,  .4 ) -- (-2.3 , - .1);
     \draw (-23mm,  3.mm) arc ( 0:180: 2.8mm);
     \draw (-2.85, .4 ) -- (-2.85,- .8);%
     \draw (-2.  ,- .35) -- (-2.   ,- .8);%
     \draw (-2.6 ,   .6) circle (4pt);
     \draw (- 2.85,  .9) node(x)  {$\heartsuit$};
     \draw (- 1.8,- .6) node(x)  {$\alpha$};
     \draw ( 4.8,- .6) node(x)  {$\beta  $};
 \draw ( 4.6, -.4 ) -- ( 4.6 , -.8);%
     \draw ( 5.45, .4) -- ( 5.45  ,- .8);%
     \draw ( 43.4mm, - 1.  mm) arc (180:360: 2.8mm);
     \draw ( 4.9,  .4 ) -- (4.9,-  .1);
     \draw ( 54.5mm, 3.0mm) arc ( 0:180: 2.8mm);
     \draw ( 4.34,1. ) -- ( 4.34 ,- .1);%
     \draw ( 5.16, 1. ) -- ( 5.16 , .6);%
     \draw ( 5.2 ,   .6) circle (4pt);
 \end{tikzpicture}\end{centering}\caption{Growing system of braids.}
 \label{growsys}\end{figure}
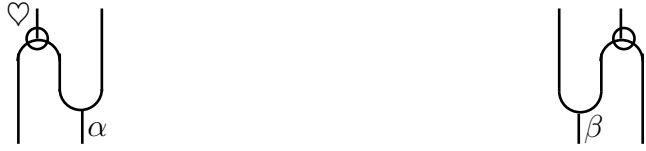

In the sequel we abbreviate for short,
\begin{gather}
(\heartsuit,\alpha)\equiv(\id\times\alpha)\circ(\heartsuit\times\id),\\
(\beta,\heartsuit)\equiv(\beta\times\id)\circ(\id\times\heartsuit).
\end{gather}

\begin{thm} Let $a$ and $b$ be unaries. Then the cooperation
$\heartsuit\equiv(a,b)$ is under $(\alpha,\beta)$-morphism iff
\begin{eqnarray*}(i)\quad&\beta\;\text{must be
$(\beta,\heartsuit)$-symmetric.}\\
(ii)\quad&\alpha\in\hom((\beta,\heartsuit),b).\\
(iii)\quad&\alpha\in\hom(\id\times a,a).\end{eqnarray*}
The cooperation $\heartsuit$ is over $(\alpha,\beta)$-morphism iff
\begin{eqnarray*}(i)\quad&\alpha\;\text{must be
$(\heartsuit,\alpha)$-symmetric.}\\
(ii)\quad&\beta\in\hom((\heartsuit,\alpha),a).\\
(iii)\quad&\beta\in\hom(b\times\id,b).\end{eqnarray*}\end{thm}

\section{Bigebra} The Boolean cogebras in combinatorics has
been considered by Joni \& Rota [1979]. However we believe that what is
most relevant to the logic is the bi-gebra as in Figure \ref{bigebra} as the
morphism, and in particular the Boolean bigebra. A general bigebra is a
triple: a pre-braid of type $\in\cat(2,2),$ binary $\in\cat(2,1)$ and cooperation $\in\cat(1,2)\simeq\cat(1,1)\times\cat(1,1).$

Therefore a bigebra has the type \begin{multline}<2,2,2,1,1>\\\simeq\cat(2,1)\times\cat(2,1)\times\cat(2,1)\times\cat(1,1)\times\cat(1,1).\end{multline}
with at least two defining relations including the Artin relation and Figure \ref{bigebra}.

\begin{figure}[!ht]
 \begin{centering}
   \begin{tikzpicture}[line width=1pt]
     \draw ( 1.3,- .0) node(x)  {$\sim$};
     \draw (- .1,   .6) -- (- .1 , .4);%
     \draw ( 0.46,  .6) -- ( 0.46 , .4);%
     \draw (- 1mm,  4.  mm) arc (180:360: 2.8mm);
     \draw (  .2,  .17) -- (  .2,- .13);
     \draw (  4.6mm, - 3.9mm) arc ( 0:180: 2.8mm);
     \draw (- .1,- .4) -- (- .1 ,- .6);%
     \draw (  .46,- .4) -- (  .46 ,- .6);
     
     \draw ( 1.91,  .2 ) -- ( 1.91,-  .1 );
     \draw ( 2.16, .6) -- ( 2.16 , .4);%
     \draw ( 19.1mm,-  .8 mm) arc (180:330: 3. mm);
     \draw ( 2.43,  .25) -- ( 2.9 ,-  .25);
     \draw ( 2.93, .32) -- ( 2.75,   .1);
     \draw ( 2.45 ,-  .26) -- ( 2.6,-  .1);
     \draw ( 24.6mm,  2.29mm) arc (30:180: 3. mm);
     \draw ( 2.2,- .4) -- ( 2.2 ,- .6);%
     
     \draw ( 3.1 , .6) -- ( 3.1  , .4);%
     \draw ( 28.9mm,- 2.0 mm) arc (210:360: 3. mm);
     \draw ( 3.45,  .2 ) -- ( 3.45,-  .1 );
     \draw ( 34.5mm,  1.19mm) arc ( 0:150: 3. mm);
     \draw ( 3.1,- .3) -- ( 3.1 ,- .6);%
   \end{tikzpicture}
  \end{centering}
  \caption{A bi-gebra with one pre-braid.}
 \label{bigebra}
\end{figure}
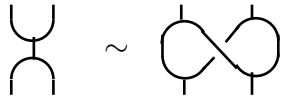

For the cloning (mitosis) and for the switch, Figure \ref{bigebra} holds for
every binary operation, \ie this bi-operation is always present
implicitely. We wish to use this relation and generalize explicitely.

\section{Dualizing object} For duality theory for Boolean algebras see [Stone 1936, Yetter 1990, Davey \& Priestley 1990, Davey 1993].

\begin{exer}[Davey 1993, p. 105, Problem 3] Which finite algebras admit
a duality?\end{exer}

Let $\Omega,A\in\text{obj\textbf{cat}}.$ An object $\Omega$ is said to
be dualizing (or schizophrenic) if for every object $A,$ there is an
isomorphism between $A$ and $\Omega^A.$ In this case an object
$\Omega^A$ is said to be the dual to $A.$ In the categories of the
finite algebras $\Omega$ is dualizing iff $|\Omega^A|=|A|.$

\begin{exer} A category of the finite sets with one modal structure.
Let $a\in A^A$ and $z\in\Omega^\Omega$ be modal (unary) structures,
\begin{gather*}(A,a),(\Omega,z)\in\text{obj\textbf{cat}},\\
f\in(\Omega,z)^{(A,a)}\quad\Longleftrightarrow\quad
f\circ a=z\circ f.\end{gather*}
Determine $\left|(\Omega,z)^{(A,a)}\right|=$

If invariant subset inv$_a\equiv\{x\in A|ax=x\}$ is not empty, and
inv$_z$ is empty then $\left|(\Omega,z)^{(A,a)}\right|=\emptyset.$

In particular let all orbits of modal structures be of cardinality 2
(then modals are unipotents and cardinality of sets must be even). If
$|\Omega|=2$ and $|A|=2n,$ then $\left|(\Omega,z)^{(A,a)}\right|=2^n.$

More complicated structure of orbits? See [Stone 1936].
\end{exer}

\begin{exer} A category of the finite sets with binaries. Let
$$\alpha\in A^{(A\times A)},\quad\omega\in\Omega^{(\Omega\times\Omega)}.
$$\end{exer}

Let $C$ be co-magma and $A$ be magma. Then $A^C$ inherit a structure of
magma with a convolution product. If magma $A$ is finite then $C^A$
inherit a structure of comagma. If $|A|<\infty$ and $\Omega$ is a
dualizing object then $A^*\equiv\Omega^A$ is comagma and this is
displayed on Figure \ref{duality}.

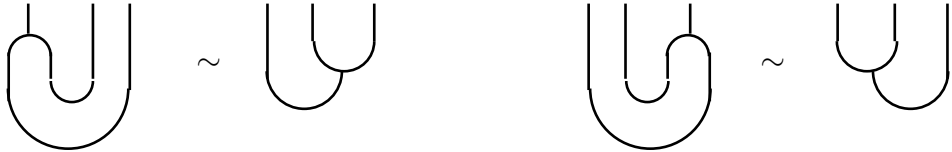
\begin{figure}[!ht]
 \begin{centering}
   \begin{tikzpicture}[line width=1pt]
     
      \draw (-4.4,1.3) -- (-4.4 , .9);%
      \draw (-3.54,1.3) -- (-3.54  ,  .3);%
      \draw (-41.mm,   2.7 mm) arc (180:360: 2.8mm);
      \draw (-4.1,  .65) -- (-4.1,   .3);
      \draw (-41mm,  6.mm) arc ( 0:180: 2.8mm);
      \draw (-4.66, .65) -- (-4.66 ,- .0);%
      \draw (-  3.05 ,1.3) -- (-  3.05  , .15);%
      \draw (-46.7mm,   1.7 mm) arc (180:360: 8mm);
      
      \draw (- 2.  ,  .5) node(x)  {$\sim$};
     
     \draw (  .2,1.3) -- ( .2,.8);%
     \draw (-1.22  ,1.3) -- (-1.22 ,.3);
     \draw (-.62,1.3) -- (- .62  ,.8);%
     \draw (-   6.  mm,  8. mm) arc ( 180:360:4 mm);
     \draw (- 12.3mm,  4.  mm) arc ( 180:360:5 mm);
          \draw ( 4.4,1.3) -- ( 4.4 , .9);%
      \draw ( 3.54,1.3) -- ( 3.54  ,  .3);%
      \draw ( 35.5mm,   2.7 mm) arc (180:360: 2.8mm);
      \draw ( 4.1,  .65) -- ( 4.1,   .3);
      \draw ( 46.5mm,  6.mm) arc ( 0:180: 2.8mm);
      \draw ( 4.66, .65) -- ( 4.66 ,- .0);%
      \draw (   3.05 ,1.3) -- (   3.05  , .15);%
      \draw ( 30.7mm,   1.7 mm) arc (180:360: 8mm);
      \draw (  5.5 ,   .5) node(x)  {$\sim$};
     \draw (6.36,1.3) -- (6.36,.8);%
     \draw (7.82  ,1.3) -- ( 7.82 ,.3);
     \draw ( 7.12,1.3) -- ( 7.12  ,.8);%
     \draw (  63.5 mm,  8. mm) arc ( 180:360:4 mm);
     \draw ( 68.3mm,  4.  mm) arc ( 180:360:5 mm);
 \end{tikzpicture}
  \end{centering}
  \caption{The product - co-product duality.}
 \label{duality}
\end{figure}

We need first to have dualizing object $\Omega=\cat(0,0).$ Maybe we must allow that
$|\Omega^A|>|A|$ ?, in order to have the perfect duality (\ie bijection)
in binary operations? That is to every binary on $A$ corresponds just
unique cobinary on $|\Omega^A|$ ?

\section{The Wigner \& Eckart problem:\\the Clebsch-Gordan coefficients}
A left action of a bi-magma or a bi-monoido $\cM$ is represented on
Figure \ref{bimagma}: the binary tree on the left is an action on a single
$\cM$-set and the next tangle represent the action of $\cM$ on
cartesian (or tensor) product of two $\cM$-sets. This last action
depends on co-product.

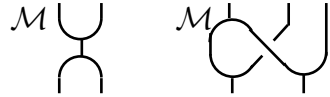
\begin{figure}[!ht]
 \begin{centering}
   \begin{tikzpicture}[line width=1pt]
     \draw (-  .5, .4) node(x)  {$\cM$};
     \draw (- .1,   .6) -- (- .1 , .4);%
     \draw ( 0.46,  .6) -- ( 0.46 , .4);%
     \draw (- 1mm,  4.  mm) arc (180:360: 2.8mm);
     \draw (  .2,  .12) -- (  .2,- .12);
     \draw (  4.6mm, - 3.9mm) arc ( 0:180: 2.8mm);
     \draw (- .1,- .4) -- (- .1 ,- .6);%
     \draw (  .46,- .4) -- (  .46 ,- .6);
     
     \draw ( 1.91,  .2 ) -- ( 1.91,-  .1 );
     \draw ( 2.16, .6) -- ( 2.16 , .4);%
     \draw ( 19.1mm,-  .8 mm) arc (180:330: 3. mm);
     \draw ( 2.43,  .25) -- ( 2.9 ,-  .23);
     \draw ( 2.94, .31) -- ( 2.75,   .1);
     \draw ( 2.45 ,-  .26) -- ( 2.6,-  .1);
     \draw ( 24.6mm,  2.29mm) arc (30:180: 3. mm);
     \draw ( 2.2,- .4) -- ( 2.2 ,- .6);%
     \draw ( 1.7,  .4) node(x)  {$\cM$};
    
     \draw ( 2.95 , .6) -- ( 2.95  , .3);%
     \draw ( 28.9mm,- 2.0 mm) arc (210:360: 3. mm);
     \draw ( 3.45,  .6 ) -- ( 3.45,-  .1 );%
     \draw ( 3.15,- .35) -- ( 3.15,- .6);%
   \end{tikzpicture}
  \end{centering}
 \caption{An action of a bi-magma $\cM$ on left $\cM$-set and a co-product dependent action of $\cM$ on the cartesian product of two left $\cM$-sets.}
 \label{bimagma}
\end{figure}

Every automorphism (a permutation) $\alpha\in\aut\,\cM$ of a magma,
give rise to the another action
$\gamma\rightarrow(\alpha\times\id_S)\circ\gamma.$ Two actions $\gamma$
and $\gamma'$ are equivalent if exists a bijection $\varphi$ such that
$\varphi\circ\gamma=\gamma'\circ(\id_{\cM}\times\varphi).$

The tangle on Figure \ref{bimagma} is the \textit{definition} of the (cartesian)
product of two $\cM$-sets resulting again in a $\cM$-set.
Therefore the product of two $\cM$-sets depends on a pre-braid and is
defined in terms of a coproduct on Figures \ref{cooper} and \ref{bimagma}.

An invariant operator $T$ (in linear algebra known as a tensor
operator) with respect to bi-magma $\cM$ is defined as a bin-ary
morphism of left $\cM$-sets and this definition is shown on Figure \ref{morinv}
. In particular $T$ is a $\cM$-invariant (commuting with the action
of $\cM$) binary operation of one-sided $\cM$-sets.

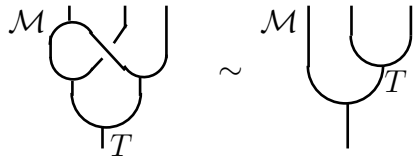
\begin{figure}[!ht]
 \begin{centering}
   \begin{tikzpicture}[line width=1pt]
     \draw (   .6, .4) node(x)  {$\cM$};
     \draw ( 1.85 ,-1.25) node(x)  {$T$};
     \draw ( 3.3 ,- .3) node(x)  {$\sim$};
    
     \draw (  .91,  .2 ) -- (  .91,-  .1 );
     \draw ( 1.16, .6) -- ( 1.16 , .4);%
     \draw (  9.1mm,-  .8 mm) arc (180:330: 3. mm);
     \draw ( 1.43,  .25) -- ( 1.9 ,-  .26);
     \draw ( 1.92, .35) -- ( 1.75,   .1);
     \draw ( 1.45 ,-  .26) -- ( 1.6,-  .1);
     \draw ( 14.6mm,  2.29mm) arc (30:180: 3. mm);
     \draw ( 1.2 ,- .4) -- ( 1.2  ,- .6);%
     
     \draw ( 1.9 , .6) -- ( 1.9  , .3);%
     \draw ( 18.9mm,- 2.0 mm) arc (210:360: 3. mm);
     \draw ( 2.45,  .6) -- ( 2.45,-  .1 );
     \draw ( 2.1,- .35) -- ( 2.1 ,- .6);%
     \draw ( 11.9mm,- 5.6 mm) arc (180:360: 4.6 mm);      
     \draw ( 1.6,-1.0 ) -- ( 1.6 ,-1.3);%
     \draw ( 3.95 ,  .4) node(x)  {$\cM$};
     \draw ( 5 .5,- .4) node(x)  {$T$};

     \draw ( 4.91, .6) -- (4.91, .2);%
     \draw ( 4.34  , .6) -- ( 4.34 ,-  .22);
     \draw (5.7 , .6) -- ( 5.7   , .2);%
     \draw (   49.  mm,  2. mm) arc ( 180:360:4 mm);
     \draw (  43.3mm,- 2.  mm) arc ( 180:360:5 mm);
     \draw ( 4.86,- .7) -- (4.86,-1.3);%
     
   \end{tikzpicture}
  \end{centering}
  \caption{The morphism $T$ of left $\cM$-sets: $T$ is
 $\cM$-invariant.}
 \label{morinv}
\end{figure}

The invariant operator $T$ is operating on (an abelian?) category of
all (one-sided) $\cM$-sets. A category of $\cM$-sets together with a
binary morphism $T$ is $\cM$-set binary-algebra. Figure \ref{morinv} has the
three other interpretations:
\begin{description}
\item[Imprimitivity] A system of imprimitivity [Weyl 1931]. This is
the case if $\cM$ is a group and if $T$ is a `projection-valued
measure'.
\item[Measuring] An action of $\cM$ is measuring bin-ary multiplication
$T.$
\item[Distribuing] An action of $\cM$ distribute $T.$
If an action of $\cM$ is denoted by dot $\cdot$ and $T=+,$ then Figure \ref{morinv} tells,$$\forall\,c\in\cM,\quad(c\cdot m)+(c\cdot n)=c\cdot(m+n).$$
\item[Evaluation] For example in $*$-autonomous category.
\end{description}

A $\cM$-invariant binary $T$ depends on a co-product
$\triangle.$ An open problem is the explicit determination of the
invariant operators $T$ for the given bimagma $\cM.$

\section{Realizations. Models in sets. Representation}
If $A,B,C$ are objects in cartesian closed category then
\be (A\times B)^C=(A^C)\times(B^C),\quad A^{(B\times C)}=
\left(A^B\right)^C.\label{ccc}\end{equation}

Throughout this paper all sets are finite.
Let $S$ be nonempty set, $s\equiv|S|\in\N.$ The multiple cartesian
products are denoted as
$$S^{\times n}\equiv
\underbrace{S\times S\times\ldots\times S}_{n\;\text{times}},$$
$$\text{with}\quad S^{\times
0}\equiv\text{one element set},\quad S^{\times 1}\equiv S,\quad
S^{\times 2}\equiv S\times S,\quad\text{etc.}$$
A clone (language, operad) of operations (graphs, plants) of type
$m\mapsto n$ can be realized (a model) inside of the power sets
\begin{equation}\text{clone}\quad\ra\quad(S^{\times n})^{(S^{\times m})}.
\label{clone}\end{equation}
Representation of a clone is a contravariant functor. If map
\eqref{clone} is injective than a set $S$ is said to be a carrier for a
clone.

An endomap of type $2\mapsto 2,$ can be realized for example as a map
$\sigma:S\times S\ra S\times S,$
$\sigma\in\End_\cF(S\times S)$ is
said to be \textit{pre-braid} if $\sigma$ solve the Artin prebraid
relation $3\mapsto 3$ represented by tangles on Figure \ref{artbraid}.

\section{The {\L}ukasiewicz logic} The three valued logic $\{$true,undefined,false$\},$ [{\L}ukasiewicz 1918, Post 1921, Kleene 1952].

\begin{table}[ht]
\caption{Some modal (unary) operations among 23 primitives}\begin{center}
\begin{tabular}{c|ccc|cc}\hline
&true&undefined&false&\bt minimum\\polynomial\et&\\\hline\\
possibility&t&t&f&$P^2=P$&{\L}ukasiewicz\\
necessity&t&f&f&$N^2=N$&{\L}ukasiewicz\\
contingent&f&t&f&$C^2=$const&\\
rotation&u&f&t&$R^3=$id&Post 1920\\
pseudo&f&f&t&$H^3=H$&\bt Heyting\\ 1930, 1966\et\\
negation&f&u&t&unipotent&Kleene 1938\\
knowledge&&&&&Kleene 1952\\
belief&&&&&Kleene 1952\\\hline\end{tabular}\end{center}\end{table}

\end{document}